\documentclass[10pt]{amsart}
\usepackage{amssymb}

\newcommand{\res}{\upharpoonright}
\newcommand{\dbullet}{\,\raisebox{-1.4pt}{\( \bullet \)}\,}
\newcommand{\dowd}{\, .\,}

\theoremstyle{plain}
\newtheorem{theorem}{Theorem}[section]

\theoremstyle{remark}

\numberwithin{equation}{section}

\begin{document}
\title[Ramsey theory: foundations and connections with dynamics]{Recent developments in finite Ramsey theory:\\ foundational aspects and connections with dynamics}

\author{S{\l}awomir Solecki}

\thanks{Research supported by NSF grant DMS-1266189.}

\address{Department of Mathematics\\
University of Illinois\\
1409 W. Green St.\\
Urbana, IL 61801, USA}

\email{ssolecki@math.uiuc.edu}

\begin{abstract} We survey some recent results in Ramsey theory. We indicate their connections with topological dynamics.
On the foundational side, we describe an abstract approach to finite Ramsey theory. We give one new application of the abstract approach through which
we make a connection with the theme of duality in Ramsey theory. We finish with some open problems.
\end{abstract}

\maketitle

\section{Ramsey theory and topological dynamics}\label{S:beg}

Recent years have seen a renewed interest in Ramsey theory that lead to advances both in proving new concrete Ramsey results and in
developing the foundational aspects of the theory. To a large extent this interest in Ramsey theory was sparked by the discovery of its close connections
with topological dynamics and especially with the notion of extreme amenability and related to it problem of computing universal minimal flows of topological groups.
A topological group is called {\bf extremely amenable} if
each continuous action of it on a compact (always assumed Hausdorff) space
has a fixed point. First such groups were discovered by Herer and Christensen \cite{HeCh} using functional analytic methods.
It was then shown by Veech \cite{Vee} that extremely amenable groups cannot be locally compact.
It turned out, however, that some very interesting groups are extremely amenable;
for example, Gromov and Milman \cite{GrMil} showed that the unitary group of a separable infinite
dimensional Hilbert space, taken with the strong operator topology and with composition as the group operation, is extremely amenable. The proof in
\cite{GrMil} of this theorem used probabilistic methods of concentration of measure through the notion of L{\'e}vy group. (L{\'e}vy groups are topological groups
possessing an increasing sequence of compact subgroups with dense union and with concentration of measure exhibited by the sequence
of the normalized Haar measures on the compact subgroups.) Concentration of measure grew to be one of the two main methods used in
proving extreme amenability.

It was not until Pestov's paper \cite{Pes} that the second general method---Ramsey theory---was discovered.
Pestov showed that the group of all increasing bijections from
$\mathbb Q$ to itself, with pointwise convergence topology and composition as the group operation,
is extremely amenable. His proof used the classical Ramsey theorem in a way that
appeared, as it turned out correctly, fundamental. Pestov's article was followed by
two papers by Glasner and Weiss \cite{GW} and \cite{GW2}, one of which \cite{GW2} used the dual Ramsey theorem of Graham and Rothschild, see
Theorem~\ref{T:GrRo} below, to determine the universal minimal flow of the group of all homeomorphisms of the Cantor set.
The full and unexpectedly tight connection between extreme amenability and Ramsey theory was then established by Kechris, Pestov, and Todorcevic in \cite{KPT}.

Theorem~\ref{T:KPT} below is the main result of the theory found in \cite{KPT}. (Paper \cite{NVT} contains some further developments.)
Recall that a topological group is {\bf non-Archimedean} if it has
a basis at the identity consisting of open subgroups. A structure, which we understand in the sense of Model Theory,
is {\bf ultrahomogeneous} if each isomorphism between
two finite substructures extends to an automorphism of the whole structure, and it is {\bf locally finite} if its finitely generated substructures are finite.
A class $\mathcal F$ of finite structures is said to have the {\bf Ramsey property} if
for any positive integer $c$ any two structures $A$ and $B$ in $\mathcal F$ there is a structure $C$ in $\mathcal F$ such that for each coloring with $c$ colors of
all substructures of $C$ isomorphic to $A$ there is a substructure $B'$ of $C$
isomorphic to $B$ such that substructures of $B'$ isomorphic to $A$ get the same color.

\begin{theorem}[Kechris--Pestov--Todorcevic \cite{KPT}]\label{T:KPT}
Let $G$ be a non-Archimedean, second countable, completely metrizable group.
Then $G$ is extremely amenable if and only if  $G$ is isomorphic to the group of all automorphisms  (taken with the pointwise convergence topology and
composition as the group operation) of
a countable, ultra-homogeneous, locally finite structure $A$ such that
\begin{enumerate}
\item[---] $A$ is linearly ordered and

\item[---] the class of all finite substructures of $A$ has the Ramsey property.
\end{enumerate}
\end{theorem}
For example, $\mathbb Q$ taken with its linear order is ultrahomogeneous, locally finite, linearly ordered, and the class of its finite substructures
consists of finite linear orders, which has the Ramsey property by the classical theorem of Ramsey \cite{Ram}. This gives back Pestov's theorem \cite{Pes}
mentioned above. More broadly, Theorem~\ref{T:KPT} related topological dynamics to Ramsey theory for finite structures, the latter having been developed by
Ne{\v s}et{\v r}il, R{\" o}dl \cite{NeRo}, \cite{NeRo2}, \cite{NeRo3}, Abramson and Harrington \cite{AbHa}, and others.

As it turned out, the connection from Theorem~\ref{T:KPT} suggested new Ramsey results. One, although not the only one, way this took place
was through comparisons
with the concentration of measure method. Given a group whose extreme amenability was proved using concentration of measure, one could
sometimes formulate a Ramsey statement that would yield the extreme amenability, and then ask if the Ramsey  statement itself held. On the other hand,
one could also ask for Ramsey statements that gave extreme amenability in situations to which concentration of measure did not apply.
We give below two examples, one on either side.

The following theorem for finite linearly ordered metric spaces was proved by Ne{\v s}et{\v r}il \cite{Nes2}. In its
statement by an {\bf order isometry} from a linearly ordered metric space $A$ to a linearly ordered metric space $B$ we understand a bijection from $A$ to $B$ that preserves the metric and the linear order.

\begin{theorem}[Ne{\v s}et{\v r}il \cite{Nes2}] Given a positive integer $c$ and
two finite linearly ordered metric spaces $A$ and $B$, there exits a finite linearly ordered metric space $C$
such that for each coloring with $c$ colors of all subspaces of $C$ order isometric to $A$ there exists a subspace $B'$ of $C$ order isometric to $B$
such that all subspaces of $B'$ order isometric to $A$ have the same color.
\end{theorem}

The theorem above implies, as shown in \cite{KPT}, that the group of all isometries of the separable Urysohn metric space taken with the pointwise convergence
topology is extremely amenable. Extreme amenability of this group
was earlier established by Pestov in \cite{Pes2} with concentration of measure methods.

To state the second theorem, also resulting from analyzing connections between Ramsey theory, concentration of measure, and extreme amenability,
consider the following notions. Let $[n]$ stand for the set $\{ 1, 2, \dots, n\}$. Given a prime number $p$,
let $({\mathbb Z}/p)^{n:l}$ be the set of all partial functions from $[n]$ to ${\mathbb Z}/p$ whose domains have at least $n - l$ elements, and let
$({\mathbb Z}/p)^n$ be the set of all functions from $[n]$ to ${\mathbb Z}/p$.
A set $L \subseteq ({\mathbb Z}/p)^{n:l}$
is called {\bf full} if there exists $h \in ({\mathbb Z}/p)^{n}$ and
$a \subseteq [n]$ with $n-l$ elements   such that for each $r\in {\mathbb Z}/p$
\[
(r + h) \upharpoonright a_r \in L
\]
for some $a\subseteq a_r \subseteq [n]$.

We now have the following Ramsey theorem. We will come back to it in the last section of the paper when discussing open problems.

\begin{theorem}[Farah--Solecki \cite{FS}]\label{T:FS}
Let $p_1, \dots,  p_k$ be prime numbers, and let $c$ be a positive integer. Then
\[
\exists l_1 \forall n_1\geq l_1 \cdots \exists l_k\forall n_k\geq l_k \hbox{ for each coloring of $\prod_{i=1}^k({\mathbb Z}/p)^{n_i:l_i}$ with $c$ colors}
\]
there exist full sets $L_1\subseteq ({\mathbb Z}/p_1)^{n_1:l_1}, \dots, L_k\subseteq ({\mathbb Z}/p_k)^{n_k:l_k}$ with
$L_1\times \cdots \times L_k$ monochromatic.
\end{theorem}

The proof of the above result uses Lovasz's method for calculating the chromatic numbers of the Kneser graphs, see \cite{Mat}. The theorem above implies that, for example, $L_0(\phi, A)$ is extremely amenable. The group
$L_0(\phi, A)$ is the completion of the group of all continuous functions (with pointwise addition) from the Cantor set $2^{\mathbb N}$
to a finite abelian group $A$
with respect to convergence in $\phi$, where $\phi$ is a diffuse {\em submeasure} on all closed-and-open subsets of
$2^{\mathbb N}$. (These groups are related to the ones considered by Herer and Christensen \cite{HeCh}.) On the other hand, it is shown in
\cite{FS} that extreme amenability of $L_0(\phi, A)$ as above cannot be proved using the concentration of measure method---such groups are not L{\'e}vy
despite possessing sequences of compact subgroups with dense unions.

There are many other examples
of recently found Ramsey theorems with application to topological dynamics; for a sample, see \cite{Jas}, \cite{NVT}, or \cite{Soc2}.

\section{Finite Ramsey theory---abstract approach}\label{S:aR}

The Kechris--Pestov--Todorcevic theory lead indirectly to rethinking of the foundations of finite Ramsey theory.
In this section, we present an abstract approach to finite Ramsey theory from \cite{Sol}.
This approach recovers most of the core Ramsey theory and makes it possible to prove new results. At the same time, it reveals
the formal algebraic structure underlying finite Ramsey theorems: there exist a single type of
structure, called Ramsey domain over a
normed composition space, that underlies Ramsey theorems. One formulates within this algebraic setting an abstract pigeonhole principle and an
abstract Ramsey statement, and proves, as the main theorem, that the pigeonhole principle implies the Ramsey statement.
This abstract Ramsey theorem, which we state at the end of this section as Theorem~\ref{T:ram},
gives particular Ramsey theorems as instances, or iterative instances, for particular Ramsey domains.

We outline the general approach in this section. For details and proofs the reader should consult \cite{Sol}. We give one new concrete application in the next
section, which will allow us to illustrate the abstract notions in a concrete situation and also to discuss the theme of duality in Ramsey theory.
We ask the reader to consult \cite{Sol}, \cite{Sol2}, and \cite{Zha} for more concrete applications.
Let us only mention here that the following theorems can be obtained as particular instances of the abstract approach to Ramsey theory,
see \cite{Sol}, \cite{Sol2}, and \cite{Zha}:
\begin{enumerate}
\item[---] the classical Ramsey theorem, see \cite{Nes};

\item[---] the van den Waerden--Hales--Jewett theorem, see \cite{Nes};

\item[---] the Graham--Rothschild theorem, \cite{GrRo71}, see also \cite{Nes};

\item[---] the versions of the two results directly above for partial rigid surjections due to Voigt, \cite{Voi}, see also \cite{Nes};

\item[---] a self-dual Ramsey theorem, \cite{Sol};

\item[---] the Milliken Ramsey theorem for finite trees, \cite{Mil}, see also \cite{Soc};

\item[---] a common generalization of Deuber's and Jasi{\'n}ski's Ramsey theorems for finite trees, \cite{Deu}, \cite{Jas};

\item[---] Spencer's generalization of the Graham--Rothschild theorem and the Ramsey
theorem for affine subspaces, \cite{Spe};

\item[---] dual Ramsey theorem for trees, \cite{Sol2}.
\end{enumerate}

\subsection{Normed composition spaces}

The algebraic structure is initially defined at the level of points and it is lifted later to the level of sets.
We describe first the point level structure.
Let $A$ and $X$ be sets. Assume we are given a {\em partial} function from $A\times A$ to $A$,
\[
(a, b)\to a\cdot b\in A,
\]
and a {\em partial} function from $A\times X$ to $X$,
\[
(a,x)\to a\dowd x\in X,
\]
such that for $a,b\in A$ and $x\in X$ if $a\dowd (b\dowd x)$ and
$(a\cdot b)\dowd x$ are both defined, then
\begin{equation}\notag
a\dowd (b\dowd x) = (a\cdot b)\dowd x.
\end{equation}
The above equation is just the usual action condition.
We assume we also have a function $\partial\colon X\to X$ and a function $|\cdot|\colon X\to L$, where $L$ is equipped with a partial order $\leq$.
The operations $\cdot$ and $\dowd$ are called a {\bf multiplication} and an {\bf action} (of $A$ on
$X$), respectively. We call $\partial$ a {\bf truncation} and $|\cdot|$ a {\bf norm}.

A structure $(A, X, \dowd, \cdot,\partial,  |\cdot|)$ as above is called a {\bf normed composition space} if the following conditions hold for $a\in A$ and $x,y\in X$:
\begin{enumerate}
\item[{\bf (i)}] if $a\dowd x$ and $a\dowd\partial x$ are defined, then
\begin{equation}\notag
\partial (a\dowd x) = a\dowd \partial x;
\end{equation}

\item[{\bf (ii)}] $|\partial x|\leq |x|$;

\item[{\bf (iii)}] if $|x|\leq |y|$ and $a\dowd y$ is defined, then $a\dowd x$ is defined and $|a\dowd x|\leq |a\dowd y|$.
\end{enumerate}

The conditions above record the interactions between pairs of objects among $\dowd$, $\partial$, and $|\cdot |$. So the action
is done by homomorphisms with respect to the truncation, by (i), the truncation does not increase the norm, by (ii), and the action respects the norm, by (iii).

We isolate one notion that will turn out to be useful later on.
Given $a,b\in A$, we say that $b$ {\bf extends} $a$ if for each $x$ for which
$a\dowd x$ is defined, $b\dowd x$ is defined as well and is equal to $a\dowd x$.

For $t\in {\mathbb N}$, we write $\partial^t$ for the $t$-th iteration of $\partial$. For a subset $P$ of $X$, we write $\partial P = \{ \partial x\colon x\in P\}$.

\subsection{Ramsey domains}

Here we lift the algebraic structure from points to subsets of $A$ and $X$.
Let ${\mathcal F}$ and ${\mathcal P}$ be families of
non-empty subsets of $A$ and $X$, respectively.
Assume we have a {\em partial} function from ${\mathcal F}\times {\mathcal F}$ to $\mathcal F$,
\[
(G,F)\to G\bullet F\in {\mathcal F},
\]
with the property that if $G\bullet F$ is defined, then it is given point-wise, that is,
$f\cdot g$ is defined for all $f\in F$ and $g\in G$ and
\[
F\bullet G = \{ f\cdot g\colon f\in F,\, g\in G\}.
\]
Assume we also have a {\em partial} function from ${\mathcal F}\times {\mathcal P}$ to
${\mathcal P}$,
\[
(F, P)\to F\dbullet P\in {\mathcal P},
\]
such that if $F\dbullet P$ is defined, then $f\dowd x$ is defined for all $f\in F$
and $x\in P$ and
\[
F\dbullet P = \{ f\dowd x\colon f\in F,\, x\in P\}.
\]

The structure $({\mathcal F}, {\mathcal P}, \dbullet, \bullet)$ as above is called a {\bf Ramsey domain} over the normed composition space
$(A,X, \dowd, \cdot, \partial, |\cdot|)$  if the following conditions hold:
\begin{enumerate}
\item[{\bf (a)}] if $F,G\in {\mathcal F}$, $P\in {\mathcal P}$, and
$F\dbullet (G\dbullet P)$ is defined, then so is $(F\bullet
G)\dbullet P$;

\item[{\bf (b)}] if $P\in {\mathcal P}$, then $\partial P\in {\mathcal
P}$;

\item[{\bf (c)}] if $F\in {\mathcal F}$, $P\in {\mathcal P}$, and
$F\dbullet\partial P$ is defined, then there is $G\in {\mathcal
F}$ such that $G\dbullet P$ is defined and for each $f\in F$
there is $g\in G$ extending $f$.
\end{enumerate}

The following two conditions on Ramsey domains are crucial in running inductive arguments.
A Ramsey domain is called {\bf vanishing} if for each $P\in {\mathcal P}$ there is $t\in {\mathbb N}$ such that the set
$\partial^t P$ has one element. It is called {\bf linear} if for each $P\in {\mathcal P}$, the set $\{ |x|\colon x\in P\}$
is a linearly ordered subset of $L$. The first one of these conditions makes it possible to start inductive arguments,
the second one is used to organize induction.

\subsection{Ramsey theorem}

Using the structure described earlier, we state here the abstract Ramsey theorem---Theorem~\ref{T:ram}.
The theorem will say that an appropriate pigeonhole principle implies
an appropriate Ramsey condition. The following statement is our Ramsey condition for a Ramsey domain $({\mathcal F}, {\mathcal P}, \dbullet, \bullet)$.
\begin{enumerate}
\item[{\bf (R)}] Given a positive integer $c$, for each $P\in {\mathcal P}$, there is an
$F\in {\mathcal F}$ such that $F\dbullet P$ is defined, and for
every coloring with $c$ colors of $F\dbullet P$ there is an $f\in F$ such that
$f\dowd P$ is monochromatic.
\end{enumerate}

For $P\subseteq X$ and $y\in X$, put
\begin{equation}\notag
P_{y} = \{ x\in P \colon \partial x = y\}.
\end{equation}
For $F\subseteq A$ and $a\in A$, let
\begin{equation}\notag
F_a = \{ f\in F\colon f\hbox{ extends }a\}.
\end{equation}

The following criterion is our pigeonhole principle, which we called local pigeonhole principle in \cite{Sol} and denoted it there by (LP).
We keep this notation here.
\begin{enumerate}
\item[{\bf (LP)}] Given a positive integer $c$, for all $P\in {\mathcal P}$ and $y\in
\partial P$, there are $F\in {\mathcal F}$ and $a\in A$ such
that $F\dbullet P$ is defined, $a\dowd y$ is defined, and for
every coloring with $c$ colors of $F_a\dowd P_{y}$ there is an $f\in F_a$ such that
$f\dowd P_{y}$ is monochromatic.
\end{enumerate}

The following is the abstract Ramsey theorem.

\begin{theorem}[Solecki \cite{Sol2}]\label{T:ram}
Let $({\mathcal F}, {\mathcal P}, \dbullet, \bullet)$ be a linear, vanishing Ramsey domain over a normed composition space.
Assume that each set in ${\mathcal P}$ is finite. Then (LP) implies (R).
\end{theorem}

\section{Duality and the dual Ramsey theorem for trees}\label{S:new}

In this section, we touch on the theme of duality.
In Ramsey theory for finite structures, Abramson--Harrington, Ne{\v s}et{\v r}il--R{\"o}dl's theorem \cite{AbHa}, \cite{NeRo3}
has a dual counterpart due to Pr{\"o}mel \cite{Pro}. This duality was made precise and shown to extend to proofs in \cite{Solstr} and \cite{Solstr2}.
In the unstructured Ramsey theory, the classical theorem of Ramsey \cite{Ram} has a dual counterpart due to Graham and Rothschild \cite{GrRo71}.
We will extend here this last instance of duality to trees and we will relate it to the concept of Galois connection.
This new concrete Ramsey result will also allow us
to give an illustration of the abstract notions presented in the previous section.

\subsection{The context for duality among trees---Galois connections}

Let $(S,\sqsubseteq_S)$ and $(T, \sqsubseteq_T)$ be two partial orders. A pair $(f,e)$ is called a {\bf Galois connection} if $f\colon T\to S$, $e\colon S\to T$, and
\begin{equation}\label{E:gin}
e\circ f \sqsubseteq_T {\rm id}_T\;\hbox{ and }\;f\circ e \sqsupseteq_S {\rm id}_S,
\end{equation}
that is, $e(f(w))\sqsubseteq_Tw$ and $v\sqsubseteq_S f(e(v))$ for all $w\in T$ and $v\in S$.
Usually the functions $e$ and $f$ in a Galois connection are assumed to be monotone. It is crucial for us, however, to use the more relaxed notion given above.
Galois connections in their abstract form were first defined by Ore in \cite{Ore}; for a comprehensive
treatment see \cite{G-S}.
As already noticed by Ore, of particular
importance are Galois connections fulfilling a strengthening of \eqref{E:gin} consisting of assuming that
equality holds in one of the two inequalities in \eqref{E:gin}.
(Ore called such connections perfect.)
Since the situation we consider, when the partial orders
are trees, is asymmetric, only one of these strengthenings is interesting---the one
with equality holding in the second formula in \eqref{E:gin}, in which case \eqref{E:gin} becomes
\begin{equation}\label{E:gin2}
e\circ f \sqsubseteq_T {\rm id}_T\;\hbox{ and }\;f\circ e = {\rm id}_S.
\end{equation}
Galois connections with \eqref{E:gin2} are sometimes called embedding--projection pairs, and
are important in denotational semantics of programming languages, see for example \cite{DrGo}.

\subsection{The notion of rigid surjection and the dual Ramsey theorem for trees}\label{Su:trd}

By a {\bf tree} we understand a finite, partially ordered set with a smallest element, called {\bf root}, and such that
the set of predecessors of each element is linearly ordered. So below, {\em all trees are non-empty and finite}.
Maximal elements of the tree order are called {\bf leaves}. We always denote the tree order on $T$ by $\sqsubseteq_T$.

Each tree $T$ carries a binary function $\wedge_T$ that assigns to each
$v,w\in T$ the largest with respect to $\sqsubseteq_T$ element $v\wedge_Tw$ of $T$ that is a predecessor of
both $v$ and $w$. By convention, we regard every node of a tree as one of its
own predecessors and as one of its own successors.

For a tree $T$ and $v\in T$, let ${\rm im}_T(v)$ be the set of all {\bf immediate successors} of $v$, and
we do not regard $v$ as one of them.
A tree $T$ is called {\bf ordered} if for each $v\in T$ we have a fixed linear order on ${\rm im}_T(v)$. Such an assignment
of linear orders defines
the lexicographic linear order $\leq_T$ on all the nodes of $T$ by stipulating that $v\leq_T w$ if $v$ is a predecessor of $w$ and,
in case $v$ is not a predecessor
of $w$ and $w$ is not a predecessor of $v$, that $v\leq_Tw$
if the predecessor of $v$ in ${\rm im}_T(v\wedge w)$ is less than or equal to the predecessor of $w$ in
${\rm im}_T(v\wedge w)$ in the given order on ${\rm im}_T(v\wedge w)$.

Let $S$ and $T$ be ordered trees. A function $e\colon S\to T$ is called a {\bf morphism} if the following
conditions hold:
\begin{enumerate}
\item[(i)] $e(v\wedge_Sw) = e(v)\wedge_T e(w)$, for all $v,w\in S$;

\item[(ii)] $e$ is monotone between $\leq_S$ and $\leq_T$, that is,
$v\leq_S w$ implies $e(v)\leq_T e(w)$, for all $v,w\in S$;

\item[(iii)] $e$ maps the root of $S$ to the root of $T$.
\end{enumerate}

Now we give the definition of functions appearing in the dual Ramsey theorem for trees.
Let $S$, $T$ be ordered trees. A function $f\colon T\to S$ is called a {\bf rigid surjection} provided
there exists a morphism $e\colon S\to T$ such that equation \eqref{E:gin2} holds. It is not difficult to see
that in this situation $f$ determines $e$ uniquely, so the definition above could be stated without invoking $e$.

Here is the dual Ramsey theorem for trees.

\begin{theorem}[Solecki \cite{Sol3}]\label{T:mn}
Let $c$ be a positive integer. Let $S,T$ be ordered trees. There exists an ordered tree $U$ such that for each coloring with $c$ colors
of all rigid surjections from $U$ to $S$ there is a rigid surjection $g_0\colon U\to T$ such that
\[
\{ f\circ g_0\colon f\colon T\to S \hbox{ a rigid surjection}\}
\]
is monochromatic.
\end{theorem}

Ramsey theorems for trees proved so far were usually stated in terms of injective morphisms $e$; see \cite{Sol2} for a survey. (An exception here is the dual Ramsey
theorem of Graham--Rothschild, which we discuss below.) Each such injective morphism $e$ is
an element of a unique pair $(f,e)$ with the pair fulfilling \eqref{E:gin2} and with $f$ being a surjective morphism. In this situation, when both
$e$ and $f$ are morphisms,
$e$ determines $f$ and $f$ determines $e$. So Ramsey theorems formulated in terms of $e$
can be equivalently stated in terms of pairs $(f,e)$ or in terms of $f$.
One could call the formulation in terms of $f$ dual.
Now, it turns out, that on the dual side, surjective morphisms $f$ are part of a much richer family of functions---rigid surjections; one
abandons the assumption that $f$ is a morphism and obtains
a Ramsey theorem for this larger class of functions. In fact, the statement for the larger class easily implies the statements for morphisms. We discuss it briefly below.

An injective morphism between ordered trees is called an {\bf embedding}. An image of a tree $S$ under
an embedding from $S$ to $T$ is called a {\bf copy} of $S$ in $T$.
The following theorem is due to Leeb, see \cite{GrRo75}.

\begin{theorem}[Leeb]\label{T:lee}
Given a positive integer $c$ and ordered trees $S$ and $T$, there is an ordered tree $U$ such that for each coloring with $c$ colors of all copies of $S$ in $U$
there is a copy $T'$ of $T$ in $U$ such that all copies of $S$ in $T'$ get the same color.
\end{theorem}

An embedding uniquely determines a copy which is the image of the embedding,
but also vice versa, a copy uniquely determines an embedding of which it is the image.
So the theorem above can be restated in terms of embeddings and
can be easily seen to be a particular case of Theorem~\ref{T:mn} by viewing an embedding $e$ as an element of pairs $(f,e)$ fulfilling \eqref{E:gin2}.

Theorem~\ref{T:mn} also generalizes the dual Ramsey theorem of Graham--Rothschild, as we indicate below.
A {\bf $k$-partition} of a set $X$ is a family of $k$ non-empty pairwise disjoint subsets of $X$ whose union is $X$.
A $k$-partition ${\mathcal P}$ is an $k$-{\bf subpartition} of an $l$-partition
${\mathcal Q}$ if each element of ${\mathcal P}$ is the union of some elements of ${\mathcal Q}$.
For $m\in {\mathbb N}$, let $[m]$ be the set $\{ 1, \dots, m\}$.
The following is the dual Ramsey theorem of Graham and Rothschild \cite{GrRo71}. (We come back to it in the last section of the paper.)

\begin{theorem}[Graham--Rothschild \cite{GrRo71}]\label{T:GrRo}
Let $c$ be a positive integer. For each $k,l$, there exists
$m$ such that for each coloring with $c$ colors of all
$k$-partitions of $[m]$ there
exists an $l$-partition ${\mathcal Q}$ of $[m]$ such that all $k$-subpartitions of ${\mathcal Q}$ get the
same color.
\end{theorem}

If $\mathcal P$ a $k$-partition of $[m]$, then we can write ${\mathcal P}= \{ p_1, \dots , p_{k}\}$ with $\min p_i <\min p_{i+1}$, for $1\leq i<k$, and
define $f_{\mathcal P}\colon [m]\to [k]$ by
\[
f_{\mathcal P}(x) = \hbox{ the unique } i\hbox{ such that } x\in p_i.
\]
Note that $[m]$, for $m\in {\mathbb N}$, is an ordered tree if we take $[m]$ with its natural order relation and
with the unique trivial
ordering of the immediate successors of each vertex. If $[m]$ is treated as a tree,
$f_{\mathcal P}\colon [m]\to [k]$ is a rigid surjection. This observation leads to a restatement of the Graham--Rothschild theorem
in terms of rigid surjections. This restatement follows easily from Theorem~\ref{T:mn} by considering ordered trees $S=[k]$ and $T=[l]$ and
viewing the resulting tree $U$ with its linear order $\leq_U$ only (and forgetting its tree order $\sqsubseteq_U$).

So Theorem~\ref{T:mn} generalizes Leeb's
result in the way Graham--Rothschild's result generalizes the classical theorem of Ramsey~\cite{Ram}, and it generalizes Graham--Rothschild's result
in the way Leeb's result generalizes Ramsey's classical theorem.

\subsection{Description of algebraic structures for the dual Ramsey theorem for trees}

We describe here concrete examples of general structures defined in Section~\ref{S:aR} that are used to prove Theorem~\ref{T:mn}.
For technical reasons, we consider only a restricted class of rigid surjections. One proves Theorem~\ref{T:mn} for this restricted class and then derives the
full version of the theorem from the restricted one. For ordered trees $S, T$, a rigid surjection $f\colon T\to S$ is called {\bf sealed} if $f^{-1}(v) = \{ w\}$, where $v$
is $\leq_S$-largest in $S$ and $w$ is $\leq_T$-largest in $T$.
For $w\in T$, let
\[
T^w = \{ v\in T\colon v\leq_T w\},
\]
and for $f\colon T\to S$ and $v\in S$, let
\[
f^v = f\res T^{e(v)},
\]
where $e\colon S\to T$ is the unique morphism with $(f,e)$ fulfilling \eqref{E:gin2}.

Now we define a normed composition space.
Let $\mathcal L$ be a family of ordered trees such that for $T\in {\mathcal L}$ and $w\in T$, we have $T^w\in {\mathcal L}$.
We will specify $\mathcal L$ later.
The sets $A$ and $X$ will be equal to each other, as will be the operations
$\dowd$ and $\cdot$. We let $A=X$ be the set of all sealed rigid surjections $g\colon T_2\to T_1$ for $T_1, T_2\in {\mathcal L}$.
Let $f, g\in A=X$. We let
$g\cdot f = g\dowd f$ be defined precisely when $f\colon T^y\to S$ and $g\colon V\to T$ for some ordered trees
$S,T, V\in {\mathcal L}$ and a vertex $y$ in $T$. We let
\begin{equation}\notag
g\cdot f = g\dowd f = f\circ g^y.
\end{equation}

For $f\in X$ whose image is a tree $S$ define $\partial f$ as follows. If $S$ consists only of its root, let
\[
\partial f = f.
\]
If $S$ has a vertex that is not a root, let $v$ be the second $\leq_S$-largest vertex in $S$, and let
\[
\partial f = f^v.
\]

Consider $\mathcal L$ as a partial order with the partial order relation given by requiring that $T_1$ be less than $T_2$ if and only if
there exists $w\in T_2$ with $T_1 = T_2^w$.
For $f\in X$,  let
\[
|f| = (\hbox{domain of }f)\in {\mathcal L}.
\]

It is easy to check that the structure $(A, X, \dowd, \cdot, \partial, |\cdot|)$ defined above is a normed composition space.

Now we define a Ramsey domain over this normed composition space. To specify $\mathcal L$, fix a family $\mathcal T$
of ordered trees such that each ordered tree has an isomorphic copy in $\mathcal T$ and such that
$T_1\cap T_2 =\emptyset$, for $T_1, T_2\in {\mathcal T}$, and let
\[
{\mathcal L} = \{ T^w\colon T\in {\mathcal T},\, w\in T \}.
\]
We consider non-empty sets $K\subseteq A=X$ for which there exist ordered trees
$T_1, T_2$ such that each element of $K$ has its domain included in $T_2$ and its image equal to $T_1$.
We require that $T_2\in {\mathcal T}$.
Since the trees in $\mathcal T$ are pairwise disjoint, each element of $K$ determines $T_2$. We define $d(K) = T_2$ and $r(K) =T_1$.
Now, let $\mathcal F$ consist of all such sets $K$ with $r(K)\in {\mathcal T}$,  and let $\mathcal P$ consist of all such
sets $K$ with $r(K)\in {\mathcal L}$. For $F_1, F_2, F\in {\mathcal F}$ and $P\in {\mathcal P}$, let
$F_1\bullet F_2$ and $F\dbullet P$ be defined precisely when $d(F_2) = r(F_1)$ and $d(P)= r(F)$, respectivly.
In these cases, we let
\[
F_1\bullet F_2 = F_1\cdot F_2\;\hbox{ and }\; F\dbullet P = F\dowd P.
\]

Again one checks that the structure defined above is a Ramsey domain that is linear and vanishing.
Condition (R) for it gives the statement of the dual Ramsey theorem for trees
(for sealed surjections); condition (LP) for it is proved using a version of the Hales--Jewett theorem, but
we will not describe this argument here. For all the proofs the reader may consult \cite{Sol3}.

\section{Further developments and problems}

We present below two groups of problems. Both of them aim at extending, in two different ways, the point of view from Section~\ref{S:aR} beyond its
original context. The first problem has to do with unifying the approach to finite Ramsey theory of \cite{Sol}, which was described in Section~\ref{S:aR},
with Todorcevic's infinite Ramsey theory of \cite{Tod}. Issues in the second group center around
proving certain analogous or finding a better understanding of Theorems~\ref{T:KPT}, \ref{T:FS} and \ref{T:GrRo} presented earlier.

First, there exists a general approach to {\em infinite} Ramsey theory given by Todorcevic in \cite{Tod} that incorporates
earlier work of Nash-Williams, Ellentuck, and Carlson, among others. Roughly speaking,
this is a theory of finding infinite sequences $(x_n)$ such that the set of all infinite sequences formed from $(x_n)$ by, for example, amalgamating or
taking subsequences or acting by a semigroup, is monochromatic. The question arises whether one can view
the approach to finite Ramsey theory outlined in Section~\ref{S:aR} as a starting point, or as the underlying layer, of the infinite Ramsey theory.
For example, given a normed composition space $(A,X, \dowd, \cdot, \partial, |\cdot|)$ as in Section~\ref{S:aR}, it is natural to consider the space of
sequences
\[
\varprojlim (X, \partial) = \{ (x_n)\in X^{\mathbb N}\colon x_n=\partial x_{n+1}\hbox{ for each } n\in {\mathbb N}\}
\]
with the induced partial action of $A$. It seems plausible that Todorcevic's theory can be recovered in spaces of the form $\varprojlim (X, \partial)$, which
would unify the two approaches.

Second, there exist certain Ramsey statements that point to a possible relationship of Ramsey theory with
combinatorial tools coming from algebraic topology as in \cite{Mat} or from fixed point theorems in convex analysis.
(A similar view is expressed by Gromov in \cite[Introduction to Section 1]{Gro}.)
We may recall that Theorem~\ref{T:FS} is proved using methods that originated with Lovasz's proof of Kneser's conjecture, which is done
with the aid of insights coming from algebraic topology around the Lefschetz fixed point theorem.
Below, we describe two other purely Ramsey theoretic statements with some intriguing additional features.
Both of them merit attention in their own right.
It would also be very interesting to see if
the combinatorial methods stemming from algebraic topology as in
\cite{Mat} can be incorporated into the approach outlined in Section~\ref{S:aR} to shed light on these or similar statements.

Moore \cite{Moo} carried out an analysis, analogous to the Kechris--Pestov--Todorcevic \cite{KPT} analysis described in Section~\ref{S:beg},
of amenability among non-Archimedean groups.
As a by-product, he uncovered a Ramsey statement relevant to amenability of well known Thompson's group $F$.
This is the group, under composition, of all piecewise linear increasing homeomorphisms of the interval $[0,1]$ whose non-differentiability points are
dyadic rationals and whose slopes are integer powers of $2$.
Moore found a Ramsey statement equivalent to amenability of $F$ (establishing which is a major problem concerning this group).
This Ramsey statement has a new feature---it involves convex combinations. We reproduce it below.

By a {\bf binary tree} we understand an ordered tree $T$, as in Section~\ref{Su:trd}, with the property that for each vertex $v\in T$
the set of its immediate successors ${\rm im}_T(v)$ has size $0$ or $2$. For $n\in {\mathbb N}$, $n>0$, let ${\mathbb T}_n$ denote
the set of all binary trees with $n$ leaves. Given a sequence of binary trees ${\vec{U}} =(U_1, \dots, U_m)$ such that the number of leaves
in all of them totals $n$ and given a tree $T$ in ${\mathbb T}_m$, let
$T({\vec{U}})$ be the tree in ${\mathbb T}_n$ that results from $T$ by attaching $U_i$ to the
$i$-th leaf of $T$, where the leaves of $T$ are numbered according to the linear order $\leq_T$ on $T$.
The root of $U_i$ is identified with the $i$-th leave of $T$ in the resulting tree. Here is the Ramsey statement formulated by Moore.

{\em For every $m$ there exists $n\geq m$ such that for each
coloring $c\colon {\mathbb T}_n \to \{ 0,1\}$ there exist non-negative numbers $\alpha_{\vec{U}}$, where $\bar U$
ranges over all $m$-tuples ${\vec{U}}= (U_1, \dots, U_m)$ of binary trees with a total of $n$ leaves, such that $\sum_{\vec{U}} \alpha_{\vec{U}} =1$ and
\[
\sum_{{\vec{U}}} \alpha_{\vec{U}}c(T(\vec{U}))
\]
is constant as $T$ varies over ${\mathbb T}_m$.}

\begin{theorem}[Moore \cite{Moo}]
The above Ramsey statement is equivalent to amenability of Thompson's group $F$.
\end{theorem}

More broadly, Moore's analysis of amenability parallel to the analysis of extreme amenability for non-Archimedean groups
lead him to a general class of Ramsey statements phrased in terms of convex combinations. At this point, no statements of this form appear to be
known that do not follow from ordinary Ramsey statements.

There is another Ramsey statement that seems to fit here.
It was formulated by Kechris, Soki{\'c} and Todorcevic \cite{KST}, and was motivated by the desire to give a Ramsey theoretic proof
of the theorem of Giordano--Pestov \cite{GP} that the group of all measure preserving transformations of the interval $[0,1]$ with Lebesgue measure,
taken with the weak topology,
is extremely amenable. The original proof in \cite{GP} used concentration of measure. What appears to be a minor modification of
the Graham--Rothschild theorem, Theorem~\ref{T:GrRo} above, yields a Ramsey statement that would imply Giordano--Pestov's result.
The statement, which was formulated in \cite{KST} and which we reproduce below, is not known to be true.

We say that a partition $\mathcal Q$ of a finite set $X$ is {\bf homogeneous} if any two sets in $\mathcal Q$ contain the same number of elements of $X$.

{\em Given a positive integer $c$, for each $k$ and $l$ there exists $m$ such that for each coloring with $c$ colors of all homogeneous
$k$-partitions of $[m]$ there exists
a homogeneous $l$-partition $\mathcal Q$ of $[m]$ such that all homogeneous $k$-subpartitions of $\mathcal Q$ get the same color.}

\end{document}